\documentclass[11pt]{article}
\usepackage[authoryear,round]{natbib}
\usepackage{graphicx}
\usepackage{latexsym}
\usepackage[left=1.5in,top=1.5in,right=1.5in,bottom=1.5in]{geometry}\usepackage{amsmath}
\usepackage{amssymb}
\usepackage{booktabs}
\usepackage{color}
\usepackage{tikz}
\usetikzlibrary{shapes}
\usetikzlibrary{arrows}
\definecolor{darkblue}{rgb}{0,0.4,0.9}
\definecolor{gray10}{rgb}{0.1,0.1,0.1}
\definecolor{gray20}{rgb}{0.2,0.2,0.2}
\definecolor{gray30}{rgb}{0.3,0.3,0.3}
\definecolor{gray40}{rgb}{0.4,0.4,0.4}
\definecolor{gray60}{rgb}{0.6,0.6,0.6}
\definecolor{gray80}{rgb}{0.8,0.8,0.8}
\definecolor{gray90}{rgb}{0.9,0.9,.9}
\definecolor{gray95}{rgb}{0.95,0.95,.95}
\definecolor{gray96}{rgb}{0.96,0.96,.96}
\definecolor{lgreen} {RGB}{180,210,100}
\definecolor{dblue}  {RGB}{20,66,129}
\definecolor{ddblue} {RGB}{11,36,69}
\definecolor{lred}   {RGB}{220,0,0}
\definecolor{nred}   {RGB}{224,0,0}
\definecolor{norange}{RGB}{230,120,20}
\definecolor{nyellow}{RGB}{255,221,0}
\definecolor{ngreen} {RGB}{98,158,31}
\definecolor{dgreen} {RGB}{78,138,21}
\definecolor{nblue}  {RGB}{28,130,185}
\definecolor{jblue}  {RGB}{20,50,100}
\definecolor{nnyellow}{RGB}{235,200,0}
\definecolor{purple}{RGB}{150, 0, 120}
\definecolor{sgGreen} {RGB}{20, 180, 50}
\definecolor{revised}{rgb}{0,0,0.9}

\newtheorem{theorem}{Theorem}

\newtheorem{lemma}{Lemma}

\newcommand{\pl}{\parallel}
\newcommand{\openr}{\hbox{${\rm I\kern-.2em R}$}}
\newcommand{\openn}{\hbox{${\rm I\kern-.2em N}$}}

\usepackage[utf8]{inputenc}
\usepackage[T1]{fontenc}

\title{Uniform Consistency of the Highly Adaptive Lasso Estimator of Infinite Dimensional Parameters}

\author{Mark J. van der Laan\thanks{This author gratefully acknowledges the support of NIH grant R01 AI074345-09} \ and Aurélien F. Bibaut \\ Division of Biostatistics, University of California, Berkeley\\ {\tt laan@berkeley.edu}
\\
} \date{\today}
\begin{document}
\maketitle
 
\begin{abstract}
Consider the case that we observe $n$ independent and identically distributed copies of a random variable with a probability distribution known to be an element of a specified statistical model. We are interested in estimating an infinite dimensional target parameter that minimizes the expectation of a specified loss function. In \cite{generally_efficient_TMLE} we defined an estimator that minimizes the empirical risk over all multivariate real valued cadlag functions with variation norm bounded by some constant $M$  in the parameter space, and selects $M$ with cross-validation. We referred to this estimator as the Highly-Adaptive-Lasso estimator due to the fact that the constrained can be formulated as a bound $M$ on the sum of the coefficients a  linear combination of a very large number of basis functions. Specifically, in the case that the target parameter is a conditional mean, then it can be implemented with the standard LASSO regression estimator. 
In \cite{generally_efficient_TMLE} we proved that the HAL-estimator is consistent w.r.t. the (quadratic) loss-based dissimilarity at a rate faster than $n^{-1/2}$ (i.e., faster than $n^{-1/4}$ w.r.t. a norm), even when the parameter space is completely nonparametric. The only assumption required for this rate is that the true parameter function has a finite variation norm. 
The loss-based dissimilarity is often equivalent with the square of an $L^2(P_0)$-type norm. In this article, we establish that under some weak continuity condition, the HAL-estimator is also uniformly consistent. 
 
\end{abstract}


\section{Introduction}
Let $O\sim P_0\in {\cal M}$ and  $\Psi:{\cal M}\rightarrow {\bf \Psi}$ be an infinite dimensional target parameter of interest, where ${\bf \Psi}=\{\Psi(P):P\in {\cal M}\}$ is the parameter space of $\Psi$. The estimand is thus given by $\psi_0=\Psi(P_0)$. We observe $n$ i.i.d. copies of $O$. We assume there exists a loss function 
$L(\psi)(O)$ such that $P_0 L(\psi_0)=\min_{\psi\in {\bf \Psi}}P_0L(\psi)$.
We assume that the loss function is uniformly bounded: 
\begin{equation}\label{boundloss}
\sup_{\psi\in {\bf \Psi}}\sup_o\mid L(\psi)(o)\mid<\infty .
\end{equation}
In the case that the loss-based dissimilarity $d_0(\psi,\psi_0)\equiv P_0L(\psi)-P_0L(\psi_0)$ is quadratic, we often also assume
\begin{equation}\label{quadraticboundloss}
\sup_{\psi\in {\bf \Psi}}\frac{P_0(L(\psi)-L(\psi_0))^2}{P_0L(\psi)-P_0L(\psi_0)}<\infty .
\end{equation}
We assume that the parameter space ${\bf \Psi}$ is a subset of $d$-variate real valued cadlag functions $D[0,\tau]$ on a cube $[0,\tau]\subset \openr^d_{\geq 0}$. A function in $D[0,\tau]$ is right-continuous with left-hand limits, and we also assume that it is left-continuous at  any point on the right-edge of $[0,\tau]$: so if $x_j=\tau_j$ for some $j\in \{1,\ldots,d\}$, then we assume that $\psi$ is continuous at such an $x$. We also assume that each function $\psi$ in the parameter space ${\bf \Psi}$ has a uniform sectional variation norm bounded by some universal $M<\infty$, but one can also select $M$ with cross-validation to avoid this assumption $\sup_{\psi\in {\bf \Psi}}\pl \psi\pl_v<\infty$ (see \cite{generally_efficient_TMLE}), in which case we only need to assume that the variation norm of each single $\psi$ is finite. We define the uniform sectional variation norm of a multivariate real valued cadlag function $\psi$ as 
\[
\pl \psi\pl_v=\psi(0)+\sum_{s\subset\{1,\ldots,d\}}\int_{0_s}^{\tau_s} \mid \psi_s(du_s)\mid,\]
where the sum is over all subsets $s$ of $\{1,\ldots,d\}$; for a given subset $s$, we define $u_s=(u_j:j\in s)$, $u_{-s}=(u_j:j\not\in s)$; and we define the section $\psi_s(u_s)\equiv \psi(u_s,0_{-s})$ that sets the components in the complement of $s$ equal to zero. 
Any cadlag function that has a bounded variation norm generates a finite measure so that integrals w.r..t this function are well defined.
We also assume that for each $\psi\in {\bf \Psi}$, $O\rightarrow L(\psi)(O)$ is a $d_1$-variate cadlag function on a compact support $[0,\tau_1]\subset \openr^{d_1}_{\geq 0}$ with universally bounded variation norm:
\begin{equation}\label{boundedvarloss}
\sup_{\psi\in {\bf \Psi}}\pl \psi\pl_v<\infty.
\end{equation}
If $O=(B,O_1)$ for a discrete variable $B\in \{1,\ldots,K\}$ and continuous component $O_1$, then one only needs to assume this for  $O_1\rightarrow L(\psi)(b,O_1)$ for each $b$.

Consider the following estimator $\hat{\Psi}:{\cal M}_{np}\rightarrow {\bf \Psi}$ defined by
\begin{equation}\label{HAL}
\hat{\Psi}(P_n)=\arg\min_{\psi\in {\bf \Psi}}P_n L(\psi).\end{equation}
In \cite{generally_efficient_TMLE} we proved that this estimator converges in loss-based dissimilarity at a rate faster than $n^{-1/2}$ to its true counterpart: 
\begin{equation}\label{rateHAL}
d_0(\psi_n,\psi_0)=P_0 L(\psi_n)-P_0L(\psi_0)=O_P(n^{-1/2-\alpha(d)}),\end{equation}
where $\alpha(d)>0$ is a specified number that behaves in the worst case as  $1/d$. The worst case corresponds with ${\bf \Psi}={\bf \Psi}_{NP}\equiv \{\psi\in D[0,\tau]:\pl \psi\pl_v<M\}$ being equal to the set of cadlag functions with variation norm bounded by $M$, while this rate will be better for smaller parameter spaces ${\bf \Psi}$ and can be expressed in terms of  the entropy of ${\bf \Psi}$.
For the case that the parameter space equals the nonparametric parameter space ${\bf \Psi}_{NP}$,  this estimator can be defined as the minimizer of the empirical risk $P_n L(\psi)$ over a linear combination of around $n 2^{d-1}$ indicator basis functions under the constrained that the sum of the absolute value of its coefficients is bounded by $M$. This is shown by using the following representation of a function $\psi\in D[0,\tau]$ with $\pl \psi\pl_v<\infty$:
\[
\psi(x)=\psi(0)+\sum_{s\subset \{1,\ldots,d\}}\int_{0_s}^{x_s} d\psi_s(u_s).\]
This representation shows that $\psi$ can be represented as an infinite linear combination of indicators $x_s\rightarrow I(u_s\leq x_s)$ indexed by a cut-off $u_s$ and subset $s$, where the sum of the absolute values of the ''coefficients'' $d\psi_s(u_s)$ equals $\pl \psi\pl_v$.
This motivated us to name it the  Highly Adaptive Lasso (HAL) estimator, and indeed in the case of a squared error or log-likelihood loss for binary outcomes  it reduces to the standard Lasso regression estimator as implemented in standard software packages, but where one runs it with a possibly enormous amount of basis functions.

For example, for the squared error loss and $\psi_0=E_{P_0}(Y\mid W)$ being a regression function, $d_0(\psi,\psi_0)=P_0(\psi-\psi_0)^2$ is the square of the $L^2(P_0)$-norm. Thus, our general convergence result will typically imply convergence in an $L^2(P_0)$ or Kullback-Leibler norm. In this article we are concerned with showing that this general HAL-estimator is also uniformly consistent under certain additional smoothness conditions. Let $\pl \psi\pl_{\infty}=\sup_{x\in [0,\tau]}\mid \psi(x)\mid $ be the supremum norm. We want to prove that 
\begin{equation}\label{uniformcons}
\pl \psi_n-\psi_0\pl_{\infty}\rightarrow_p 0. \end{equation}

\paragraph{The case that the observed data has a discrete and continuous component}
Before we proceed we demonstrate how one can apply our results to a setting in which $\psi_0$ is a function of a purely discrete component $B$ and continuous component. 
Suppose that $O=(B,O_1)$, where $B$ is discrete with finite number of values $\{1,\ldots,K\}$, and $\psi=(\psi_b:b=1,\ldots,K)$, where the components $\psi_b$ are variation independent so that ${\bf \Psi}=\prod_{b=1}^K{\bf \Psi}_b$ with ${\bf \Psi}_b$ being the parameter space of $\Psi_b:{\cal M}\rightarrow {\bf \Psi}_B$. One now assumes that for each $b$ ${\bf \Psi}_b$ is a subset of $d_b$-dimensional cadlag functions with variation norm smaller than some $M_b<\infty$.  We have $d_0(\psi_n,\psi_0)=\sum_{b=1}^K \int \{L(\psi_n)(b,o_1)-L(\psi_0)(b,o_1)\}dP_0(b,o_1)$. Suppose that $L(\psi)(b,o_1)$ only depends on $\psi$ through a $\psi_b$ and suppose that $\psi=(\psi_b:b=1,\ldots,K)$ is a variation independent parameterization. Then, $\psi_{0,b}$ is the minimizer of $\psi\rightarrow P_0L_b(\psi)$ where $L_b(\psi)(O_1)=I(B=b)L(\psi)(b,O_1)$, and $\psi_{n,b}=\arg\min_{{\bf \Psi}_b}P_0L_b(\psi)$. In addition, $d_0(\psi_n,\psi_0)=\sum_b d_{0,b}(\psi_{n,b},\psi_{0,b})$, where
$d_{0,b}(\psi_b,\psi_{0,b})=P_{0}L_b(\psi_b)-P_0L_b(\psi_{0,b})$. Thus the estimator $\psi_n$ above can then be analyzed separately as an estimator $\psi_{n,b}$ for $\psi_{0,b}$ for each $b$. 
In particular, the rate of convergence result above now applies to each $\psi_{n,b}$ with dimension $d$ replaced by $d_b$ and loss function $L_b(\psi)$. Our goal is then reduced to establishing that $\psi_{n,b}-\psi_{0,b}$ converges uniformly to zero in probability. In the sequel we suppress this index $b$, but the reader needs to know that in such applications we simply apply our results  to $\psi_{0,b}$ and $\psi_{n,b}$ with loss function $L_b(\psi_b)$, for each $b$ separately.  In order to establish our uniform consistency result, we will assume that each $P_0(B=b,\cdot)$ is a continuous measure for $O_1$, which corresponds with the stated assumption $\textbf{A2}$ below that $P_0$ is continuous on the support of $L_b$.

To establish the uniform consistency we will make the following assumptions:
\begin{description}
\item[A0]: $d_0(\psi_n,\psi_0)=o_P(1)$ and the loss function is uniformly bounded (\ref{boundloss}).
\item[A1]: $d_0(\psi,\psi_0)=0$ implies $\pl \psi-\psi_0\pl_{P_0}=0$. 
\item[A2]: $\psi_0$ is continuous on $[0,\tau]$, and $P_0$ is continuous measure on the set  of $o$-values for which $\sup_{\psi}\mid L(\psi)(o)\mid>0$.
\item[A3]: If $\psi_n$ converges pointwise to $\psi_{\infty}\in {\bf \Psi}$ on $[0,\tau]$ at each continuity point of $\psi_{\infty}\in {\bf \Psi}$, then $L(\psi_n)$ converges pointwise to $L(\psi_{\infty})$ on a support of $P_0$.  
\end{description}
 Regarding assumption $\textbf{A0}$, above  we provided sufficient assumptions that even guarantee $d_0(\psi_n,\psi_0)=O_P(n^{-1/2-\alpha(d)})$, which could thus easily be weakened, as long as we keep assuming that the loss function is uniformly bounded.
Assumption $\textbf{A1}$ is a very weak assumption. Regarding assumption $\textbf{A3}$, since $P_0$ is continuous by $\textbf{A2}$, one only needs to show that $L(\psi_n)$ converges to $L(\psi_{\infty})$ on a set that can exclude any finite or countable set. Since the number of discontinuity points of $\psi_{\infty}$ is finite or countable, the lack of convergence of $\psi_n$ at these points should not be  an issue.

We have the following theorem. 
\begin{theorem}\label{thuniformconsistent}
Let $\psi_n$ be the HAL-estimator defined by (\ref{HAL}). 
Assume $\textbf{A0}$, $\textbf{A1}$, $\textbf{A2}$ and $\textbf{A3}$. Then, $\sup_{x\in [0,\tau]}\mid\psi_n(x)-\psi_0(x)\mid \rightarrow0$ in probability as $n\rightarrow\infty$.  
\end{theorem}

 \section{Proof of Theorem \ref{thuniformconsistent}}
 Using that $\sup_{\psi\in {\bf \Psi}}\sup_o\mid L(\psi)(o)\mid<\infty$, the dominated convergence theorem combined with $\textbf{A3}$ proves the following lemma.
\begin{lemma}\label{lemmapointwise}
Assume $\textbf{A0}$ and $\textbf{A3}$.
If $\psi_n$ converges pointwise to $\psi_{\infty}\in {\bf \Psi}$ on $[0,\tau]$ at each continuity point of $\psi_{\infty}$, then $P_0L(\psi_n)-P_0L(\psi_{\infty})\rightarrow 0$.
\end{lemma}
The following lemma proves that if $d_0(\psi,\psi_0)=0$, then $\psi$ equals $\psi_0$ pointwise as well.
\begin{lemma} \label{lemmaequalsupnorm} Assume $\textbf{A1}$. 
If $d_0(\psi,\psi_0)=0$ for a $\psi,\psi_0\in D[0,\tau]$, then $\pl \psi-\psi_0\pl_{\infty}=0$. 
\end{lemma}

\medskip

{\bf Proof:} Assume $d_0(\psi,\psi_0)=0$. Suppose that $\psi-\psi_0>0$ (same for $<0$) at a point $x\in [0,\tau)$, then it will also be larger than $0$ at a small neighborhood $[x,x+\delta)$ for some $\delta>0$ due to the right-continuity of $\psi-\psi_0$. As a consequence, if $\psi-\psi_0>0$ at a point $x$, then $\pl \psi-\psi_0\pl_{P_0}>0$. By assumption $\textbf{A1}$ this implies that $d_0(\psi,\psi_0)>0$. Finally, if $x\in [0,\tau)^c\subset [0,\tau]$, then we assumed that $\psi,\psi_0$ are left-continuous, so that the same argument applies if we assume that $\psi-\psi_0>0$ at an $x$ on the right-edge of $[0,\tau]$.
This proves that $d_0(\psi,\psi_0)=0$ implies that $\psi-\psi_0=0$ on $[0,\tau]$.  
$\Box$

The following lemma establishes that our parameter space ${\bf \psi}$ is weakly compact so that each sequence has a weakly converging (i.e., poinwise) subsequence. 
In addition, if we also assume that the sequence is consistent for $\psi_0$, then the limit of this weakly converging subsequence has to equal $\psi_0$ as well.
\begin{lemma}\label{lemmapointwiseconvergence}
Assume $\textbf{A0}$, $\textbf{A1}$, $\textbf{A2}$, and $\textbf{A3}$. Any   sequence $(\psi_n:n=1,\ldots)$ in ${\bf \Psi}$ has a subsequence  $(\psi_{n(k)}:k=1,\ldots)$ so that there exists a $\psi_{\infty}\in {\bf \Psi}$ and $\psi_{n(k)}$ converges pointwise to $\psi_{\infty}$ at each continuity point of $\psi_{\infty}$.

If we also know that $d_0(\psi_n,\psi_0)\rightarrow 0$, then we have that $\pl \psi_{\infty}-\psi_0\pl_{\infty}=0$. 
\end{lemma}

\medskip

{\bf Proof:}
By \cite{hildebrandt1967} (see also lemma 1.2 in \cite{marks_phdthesis}), any cadlag function of bounded variation can be represented as a difference of two monotone cadlag functions generating positive finite measures, i.e. the analogue of cumulative distributions functions  but not bounded by $[0,1]$. Thus $\psi_n=F_n-G_n$ for monotone increasing functions $F_n,G_n\in D[0,\tau]$. 
Any sequence $(F_n:n)$ of cumulative distribution functions has a subsequence that converges weakly to a limit $F_{\infty}$, and similarly, any sequence $(G_n:n)$ has a subsequence that converges weakly to a limit $G_{\infty}$, where weak convergence is equivalent with pointwise convergence at each continuity point of the limit.
This shows that we can find a subsequence $(F_{n(k)}-G_{n(k)}:k)$ of $(F_n-G_n:n)$ and limit $\psi_{\infty}=F_{\infty}-G_{\infty}$ so that $F_{n(k)}-G_{n(k)}$ converges pointwise to $F_{\infty}-G_{\infty}$ at each point in which both $F_{\infty}$ and $G_{\infty}$ are continuous. We now want to show that the points at which $\psi_{\infty}$ are continuous are equal to the point at which both $F_{\infty}$ and $G_{\infty}$ are continuous. By the Hahn decomposition theorem both $F_{\infty}$ and $G_{\infty}$ are the sum of a continuous measure and purely discrete measure. The continuous measure corresponds with a continuous function. The discrete support of $F_{\infty}$ and $G_{\infty}$ has to be disjoint since if a measure assigns at a point both a negative and positive mass then we can replace that by just assigning a single mass that is either positive or negative.  
Thus we have shown that $(\psi_{n(k)}:k)$ converges pointwise  to $\psi_{\infty}$ at each continuity point of $\psi_{\infty}$. 

Consider now the second statement in the lemma. Suppose now that we also know that $d_0(\psi_n,\psi_0)\rightarrow_p 0$. Then we also have $d_0(\psi_{n(k)},\psi_0)\rightarrow 0$. 
By Lemma \ref{lemmapointwise}, the fact that $\psi_{n(k)}$ converges pointwise to $\psi_{\infty}$ at each continuity point of $\psi_{\infty}$ implies that $P_0L(\psi_{n(k)})-P_0L(\psi_{\infty})\rightarrow 0$. Now use that $d_0(\psi_{n(k)},\psi_0)=P_0L(\psi_{n(k)})-P_0L(\psi_{\infty})+d_0(\psi_{\infty},\psi_0)$. Since the left-hand side converge to zero, and the first term on the right-hand side converges to zero as well, this implies that $d_0(\psi_{\infty},\psi_0)=0$. By Lemma \ref{lemmaequalsupnorm}, this implies that $\pl \psi_{\infty}-\psi_0\pl_{\infty}=0$. This completes the proof of the lemma.
$\Box$

Consider our HAL-estimator $\psi_n$. Given $d_0(\psi_n,psi_0)\rightarrow_p 0$, Lemma \ref{lemmapointwiseconvergence} proves that $\psi_n$ converges pointwise to $\psi_0$ at each point in $[0,\tau]$, where $\psi_0$ is continuous. 
Thus, we have translated the consistency of $\psi_n$ w.r.t. loss-based dissimilarity into pointwise convergence.  
\begin{lemma}\label{thpointwisehal}
Let $\psi_n$ be the HAL-estimator defined by (\ref{HAL}). 
Assume $\textbf{A0}$, $\textbf{A1}$, $\textbf{A2}$ and $\textbf{A3}$. Then, $\psi_n(x)-\psi_0(x)\rightarrow_p 0$ at each 
$x\in [0,\tau]$. 
More generally, we have $\psi_n=F_n-G_n$ for $F_n,G_n$ that generate positive uniformly finite measures, $\psi_0=F_0-G_0$ for $F_0,G_0$ that generates finite positive measures, and $F_n(x)-F_0(x)\rightarrow_p 0$ and $G_n(x)-G_0(x)\rightarrow 0$ for each $x\in [0,\tau]$.  
\end{lemma}

%https://math.stackexchange.com/questions/467976/uniform-convergence-and-weak-convergence
So we have shown $\psi_n=F_n-G_n$, $\psi_0=F_0-G_0$, where $F_n,G_n$ converge pointwise to $F_0,G_0$ at each point in $[0,\tau]$.
Finally, we establish that the pointwise convergence of $F_n$ ($G_n$) to a continuous $F_0$ ($G_0$)  implies uniform convergence, thereby showing that $\psi_n$ converges uniformly to $\psi_0$ as well. 
\begin{lemma}\label{lemmapointwisetounif}
If $F_n$ is a sequence of cadlag functions that generate a positive measure on $[0,\tau]$, $F_n(x)\rightarrow F_0(x)$ for each $x\in [0,\tau]$, and $F_0$ is continuous on $[0,\tau]$, then
$\pl F_n-F_0\pl_{\infty}\rightarrow 0$.
\end{lemma}

\medskip

{\bf Proof:}
Let $\epsilon > 0$.

By Heine's theorem, since $F_0$ is continuous on the compact set $[0, \tau]$, it is uniformly continuous on $[0, \tau]$. 

By uniform continuity of $F_0$, there exists $\eta > 0$ such that for any $x, y \in [0, \tau]$, $\|x - y\| < \eta$ implies $|F_0(x) - F_0(y)| < \epsilon$. Consider a grid on $[0, \tau]$ with grid points $x_{\textbf{i}} \equiv (i_1 \eta, ..., i_d \eta )$.

Consider an arbitrary $x \in [0, \tau]$. For a certain $\textbf{i} \in \mathbb{N}^d$, $x$ falls in the hypercube $[x_\textbf{i}, x_\textbf{i + 1}]$, where $\textbf{1} \equiv (1,...,1)$.

Since $F_0$ and $F_n$ generate positive measures,
\begin{gather}\label{brackets}
F_n(x_\textbf{i}) - F_0(x_\textbf{i+1}) \leq F_n(x) - F_0(x) \leq F_n(x_\textbf{i+1}) - F_0(x_\textbf{i}).
\end{gather}

Observe that
\begin{equation}\label{upper_bound}
F_n(x_\textbf{i+1}) - F_0(x_\textbf{i}) = (F_n(x_\textbf{i+1}) - F_0(x_\textbf{i+1})) + F_0(x_\textbf{i+1}) - F_0(x_\textbf{i+1})).
\end{equation}

Since $F_n(x_\textbf{j}) - F_0(x_\textbf{j})$ converges to zero for all the $x_\textbf{j}$'s in $[0, \tau]$, and since there are a finite number of such $x_\textbf{j}$'s, there exists $n_0 > 0$ such that for all $x_\textbf{j} \in [0, \tau]$, $n > n_0$, $|F_n(x_\textbf{j}) - F_0(x_\textbf{j})| < \frac{\epsilon}{2}$.

Therefore, going back to \eqref{upper_bound} and using this latter fact and the uniform continuity, we have $F_n(x_\textbf{i+1}) - F_0(x_\textbf{i}) \leq \frac{\epsilon}{2}$ for any $n > n_0$.

Since we can apply the exact same arguments to the lower bound in \eqref{brackets}, we have that for $n > n_0$, 
\begin{equation}
-\frac{\epsilon}{2} \leq F_n(x) - F_0(x) \leq \frac{\epsilon}{2}.
\end{equation}

Since $n_0$ does not depend on $x$, we have proved uniform convergence of $F_n$ to $F_0$ over $[0, \tau]$.
$\Box$

\bibliography{biblio}

\begin{thebibliography}{3}
\providecommand{\natexlab}[1]{#1}
\providecommand{\url}[1]{\texttt{#1}}
\expandafter\ifx\csname urlstyle\endcsname\relax
  \providecommand{\doi}[1]{doi: #1}\else
  \providecommand{\doi}{doi: \begingroup \urlstyle{rm}\Url}\fi

\bibitem[Hildebrandt(1963)]{hildebrandt1967}
Theophil~Henry Hildebrandt.
\newblock \emph{Introduction to the theory of integration}.
\newblock Academic Press, 1963.

\bibitem[van~der Laan(1993)]{marks_phdthesis}
Mark~J. van~der Laan.
\newblock \emph{Efficient and Inefficient Estimation in Semiparametric Models}.
\newblock PhD thesis, 1993.

\bibitem[van~der Laan(2015)]{generally_efficient_TMLE}
Mark~J. van~der Laan.
\newblock A generally efficient targeted minimum loss based estimator.
\newblock \emph{U.C. Berkeley Division of Biostatistics Working Paper Series},
  2015.

\end{thebibliography}

 \end{document}